\documentclass[11pt,a4paper]{article} %\documentclass[11pt,a4paper,leqno]{article}

\usepackage[latin1]{inputenc}
 \usepackage{verbatim}
 \usepackage{enumerate}
  \usepackage{array} \usepackage{amsmath}
   \usepackage{amssymb} \usepackage{amsfonts}
   \usepackage{theorem} \usepackage{latexsym}
    \usepackage{color}
 \newcommand{\ve}{\varepsilon}
  \newcommand\NN{\mathbb{N}}
   
   \newcommand\CC{\mathbb{C}}
\newcommand{\no}{\noindent} \newcommand{\beq}{\begin{equation}} \newcommand{\eeq}{\end{equation}}

\begin{document}

\begin{center}{\bf On overconvergent subsequences of closed to rows classical Pad{\'e} approximants}\end{center}
\begin{center}{\sl Ralitza K.Kovacheva \\Institute of Mathematics and Informatics,  Bulgarian Academy of Sciences,
 Acad. Bonchev str. 8, 1113 Sofia, Bulgaria, rkovach@math.bas.bg}\end{center}

\bigskip \no {{\bf Abstract:} {\it Let $f(z) := \sum f_\nu z^\nu$ be a power series with positive radius of convergence. In the present paper, we study the phenomenon of overconvergence of sequences of classical Pad{\'e} approximants $\{\pi_{n,m_n}\}$
 associated with $f,$ where $ m_n\leq m_{ n+1}\leq m_n+1$
 and $m_n = o(n/\log n),$ resp. $m_n = 0(n)$ as $n\to\infty.$
 We extend classical results by J. Hadamard and A. A. Ostrowski  related to
  overconvergent Taylor polynomials, as well as results by G.
   L{\'o}pez Lagomasino and A. Fern{\'a}ndes Infante  concerning
    overconvergent  subsequences of  a fixed  row of the Pad{\'e} table.}

\medskip \no{{\bf MSC:} 41A21,  41A25, 30B30}

\no{{\bf Key words:} {\it Pad{\'e} approximants, overconvergence,
 meromorphic continuation,  convergence in {$\sigma$-content.} }

\bigskip \no{\bf Introduction}

\no   Let \begin{equation}f(z):= \sum_{j=0}^{\infty}{f_j}{z^j}\label{sum}\end{equation}
  be a  power series with positive radius of convergence
 $ R_0(f):= R_0, \,R_0 > 0$. By $f$ we will denote not only
 the sum of $f$ in $D_{R_0}:= \{z, |z| < R_0\}$ but also the holomorphic (analytic and single valued) function determined
 by the element $(f, D_{R_0}).$
Fix a nonnegative integer $m\, (m\in\NN)$ and denote by $R_m(f):=R_m$ the {\it radius of $m-$meromorphy } of  $f$:  that is the radius of the largest disk centered at the zero into which  the power series $f$ admits a continuation as a meromorphic function with no more than $m$ poles (counted with regard to their multiplicities). As it is known (see \cite{gonchar1}), $R_m > 0$ iff $R_0 > 0.$ Analogously, we define {\it the radius of meromorphy} $R(f)$ as   the radius of the greatest disk $D_R$  into which $f$ can be extended as a meromorphic function in $\CC.$
 Apparently, $R(f)\geq R_m \geq R_0$.   We denote the meromorphic continuations again by $f.$

Given  a pair  $(n,m),\, n,m\in\NN$, let $\pi_{n,m}$ be the classical Pad{\'e} approximant of $f$ of order $(n,m)$. Recall that (see \cite{Pe}) ${\pi}_{n,m} = p/q$, where $p,q$ are polynomials of degree $\leq n,m$ respectively and satisfy $$(fq-p)(z) = 0(z^{n+m+1}).$$ As it is well known (\cite{Pe}), the Pad{\'e}
 approximant $\pi_{n,m}$ always exists and is uniquely determined by the conditions above.

Set $$\pi_{n,m}:= P_{n,m} /Q_{n,m},$$ where  $P_{n,m}$ and $ Q_{n,m}$ are relatively prime polynomials (we write $(P_{n,m}, Q_{n,m}) = 1).$

We recall the concept of {\it convergence in $\sigma$-content }
 (cf. \cite{Go3}).
 Given a set  $e\subset \CC,$ we put
    \begin{gather*}
    \sigma(e) := {\rm inf} \left\{\sum_\nu |V_\nu| \right\}
    \end{gather*}
where the infimum is taken over all coverings $\{\bigcup V_\nu\}$ of $e$ by disks  and $|V_\nu|$ is the diameter of the disk $V_\nu$.
 Let $\Omega$ be an open set in $\CC$ and $\varphi$ a function defined in
$\Omega$ with values in $\overline{\CC}$. The sequence of functions $\{\varphi_n\}$, rational  in $\Omega$, is said  to converge in {\it $\sigma-$content to a function $\phi$ inside $\Omega,$} if for each compact set $K\subset\Omega$ and $\ve > 0$ $\sigma\{z, z\in K, |\varphi_n(z) - \varphi(z)| > \ve \} \to 0$ as $n\to 0.$ The sequence  $\{\varphi_n\}$ converges to $\varphi,$ as $n\to\infty,$  {\it $\sigma-$almost uniformly inside $\Omega$,} if for  any compact set  $K \subset \Omega$ and every $\varepsilon > 0$ $\{\varphi_n\}$ converges to $\varphi$ uniformly in the $max-$ norm on a set of the form $K\setminus K_\ve,$ where $\sigma(K_\varepsilon) < \varepsilon$.  Analogously, we define {\it convergence in Green's capacity } and  {\it convergence almost uniformly in Green's capacity} inside $\Omega$. It follows from Cartan's  inequality $\hbox{cap}(e)\geq C \sigma (e)$
   (see \cite{landkoff}, Chp.3) that convergence in capacity implies
   $\sigma-$convergence. The reader is referred for details to \cite{Go3}.

The next result may be found in \cite{gonchar1}.

\smallskip \no{\bf Theorem 1, (\cite{gonchar1}):} {\it Given a power series (1) and a fixed integer $m\in\NN,$ suppose that $0 < R_m< \infty.$

Then the  sequence $\{\pi_{n,m}\},\, n\to\infty, m-$fixed converges $\sigma-$almost uniformly to $f$ inside $D_{R_m}$ and
 $$\limsup_{n\to\infty} \Vert f - \pi_{n,m}\Vert_{K(\ve)}^{1/n} = \frac{\max_K |z|}{R_m}$$
 for any compact subset $K$ of $D_{R_m}$ and $\ve > 0.$ }

\no(here $\Vert...\Vert_K$ stands for the $\max-$ norm on $K$.)

Theorem 1 generalizes the classical result of Montessus de Ballore about rows in the Pad{\'e} table (\cite{montessus}).

\smallskip { In the present paper,  we will be concentrating on
 the case $\limsup_{n\to\infty}m_n = \infty$. If the sequence
  $\{m_n\}$ increases "slowly enough", i,e, if  $m_n = o(n),$ (resp. $m_n = o(n/\log n) $
  as $n\to\infty,$) then the following result is valid:

\smallskip

 \no{\bf Theorem 2, (\cite{Go3}, Chpt.3):}
 {\it Given $f$ with $0 < R(f) < \infty,$ let $m_n = o(n/log n),\, n\to\infty$.
 Then the sequence  $\{\pi_{n,m_n}\}$ converges
 $\sigma-$almost uniformly to $f$ inside $D_{R(f)}.$

 In case  $m_n = o(n), \, n\to\infty$, the sequence $\{\pi_{n,m_n}\}$
 converges to $f$  in Green's capacity inside $D_{R(f)}.$

For any compact set $K\in D_{R(f)}$ and any $\ve > 0$
 $$\limsup_{n\to\infty} \Vert f - \pi_{n,m_n}\Vert_{K(\ve)}^{1/n} \leq  \frac{\max_K |z|}{R(f)}$$ }

In \cite{blkov}, the question about specifying the speed of convergence above was posed. It was shown that for a class of functions the following result is valid:

\smallskip \no{\bf Theorem 3, (\cite{blkov}):} {\it Given $f$ with $0 < R(f) < \infty,$ let $ m_n\leq n, m_n\leq m_{n+1}\leq m_n+1,\, m_n=o(n/\log n),\, n\to\infty. $
 Suppose that $f$ has a multivalued singularity on $\partial{D_{R(f)}}.$

 Then the sequence $\{\pi_{n,m_n}\}$  converges $\sigma-$almost uniformly to $f$
  inside  the disk $D_{R(f)}$ and \beq\limsup_{n\to\infty} \Vert f -
 \pi_{n,m_n}\Vert_{K(\ve)}^{1/n} =\frac{\max_{z\in K} |z|}{R(f)}\label{blattkov} \eeq for every compact set $K\subset D_{R(f)}$ and every $\ve
 >0$.}}

{ Research devoted to   imposing  weaker conditions on  the growth of  the sequence $\{m_n\}$  as $n\to\infty$ was carried out by  H.P. Blatt.
  It follows from his results that the statement of Theorem 2 remains valid if $m_n = o(n)$ as $n\to\infty.$ Furthermore,  the sequence
  $\{\pi_{n,m_n}\}$ converges almost uniformly to $f$ in capacity inside $D_{R(f)}$ (
   see  the comprehensive paper  \cite{blatt}). %\end{document}

 Let now the sequence $\{m_n\}$ of positive integers satisfy
 the conditions $m_n\leq n,\, m_n\leq m_{n+1}\leq m_n+1$.
  Set $$\pi_{n,m_n} = \frac{P_{n,m_n}}{Q_{n,m_n}}:=\pi_n = P_n/Q_n,$$
 {where} $(P_{n},Q_{n}) = 1;$   $\hbox{deg}P_n\leq n,\,\hbox{deg}Q_n \leq m_n.$

 Denote by $\tau_{n,m_n}:=\tau_n$ {\it the  defect}
 of $\pi_n;$ that is $\hbox{min}(n-\hbox{deg}P_n, m_n-\hbox{deg}Q_n). $
  Then the order of the zero of $f(z)-\pi_{n}$ at $z = 0$ is not less than $n+m_n+1-\tau_{n}$ (see \cite{Pe}); in other words
 \beq f(z)-\pi_{n}(z) = 0(z^{n+m_n+1-\tau_{n}})\label{0001}.\eeq
 Following the terminology of G. A. Baker, Jr.  and P. Gr. Morris (see \cite{baker}, p. 31),
 we say that the rational function $\pi_{n}$ exists iff $\tau_{n} = 0.$

The zeros $\zeta_{n,l}, 0 \leq l\leq m_n$ of the polynomial $Q_n$ are called {\it free poles} of the rational function $\pi_{n}.$ Let $\mu_n$ be the exact degree of $Q_n,\, \mu_n\leq m_n.$ We shall always normalize $Q_{n}$ by the condition \begin{equation} Q_{n}(z) =
 \prod(z-\zeta_{n,l}^*)\prod(1 - \frac{z}{\tilde\zeta_{n,l}}) \label{a1}\end{equation}
where $|\zeta_{n,l}^*| < 2R(f)$ and $|\tilde\zeta_{n,l}| \geq 2R(f).$

Set
 \begin{equation}P_{n}(z) = a_n z^{\hbox{deg} P_n} +\cdots.\label{a2}\end{equation}

{Suppose that $\tau_n > 0$ for some $n\in\NN$ (comp. (\ref{0001})). Then, by the block structure of the Pad{\'e} table }(see \cite{Pe})
 {   $\pi_{n-l,m-k}\equiv \pi_{n,m}$
 if  $\max(k, l) \leq \tau_{n}$. Suppose that  $f(z)-\pi_n(z) = B_nz^{n+m_n+1-\tau_n}$ with $B_n\not= 0$.
 Then $\tau_{n+1} = 0$ and
 $\pi_n\not= \pi_{n+1}.$ }}
 The definition of Pad{\'e} approximants leads to
\begin{equation}\pi_{n+1}(z) - \pi_{n}(z) = A_{n}\frac{z^{n+m_n+1-\tau_n}}{Q_n(z)Q_{n+1}(z)}, \label{b2}\end{equation}
 where
\beq A_{n} = \left\{\begin{array}{ll} a_{n+1}(\prod\frac{-1}{\tilde\zeta_{n, k}})-a_{n}(\prod\frac{-1}{\tilde\zeta_{n+1}}),& m_{n+1}=m_n+1\\ a_{n+1}(\prod\frac{-1}{\tilde\zeta_{n, k}}),& m_{n+1} = m_n\\ \end{array} \right. \label{c1}\eeq

  It was shown in \cite{gonchar1}, Eq. 33 (see also \cite{{vavprsue}})
 that for a fixed $m\in \NN$ the Pad{\'e} approximant
 $\pi_{n,m}, m-$ fixed converges, as $n\to\infty$,  together with the series
 $\sum_{n=1}^{\infty}\frac{A_nz^{n+m+1-\tau_n}}{Q_n(z)Q_{n+1}(z)}$,
 i.e.,
 $$f(z) - \pi_{n,m}(z) = \sum_{k=n}^{\infty}\frac{A_kz^{k+m+1-\tau_k}}{Q_k(z)Q_{k+1}(z)},$$
  where  $\limsup|A_n|^{1/n} = 1/R_m$.

 It is easy to check that
  under the conditions of Theorem 3
  an analogous  result holds also
 for sequences $\{\pi_{n,m_n}\},\, \{m_n\} -$ as in Theorem 3 (compare with  (\ref{b22}) below).
  In other words,
\begin{equation}f(z) - \pi_n(z) = \sum_{k=n}^\infty \frac{A_kz^{k+m_k+1-\tau_k}}{Q_k(z)Q_{k+1}(z)}\label{series}\end{equation} and $$\limsup_{n\to\infty} |A_n|^{1/n} = 1/R(f).$$ It follows from (\ref{c1}) that, under the above conditions on the growth of the sequence $\{m_n\}$  as $n\to\infty$ $$\limsup_{n\to\infty} |a_n|^{1/n} =  1/R(f)$$

Basing  on the block structure of the Pad{\'e} table (see \cite{Pe}), we will be assuming throughout the paper,    that
 $\tau_{n} = 0$ for all $n\in\NN.$ Also, for the sake of simplicity,
 we assume  that $\hbox{deg} P_n = n$ for all $n\in\NN$.

\smallskip

{ Let $f(z)=\sum_{n=0}^\infty f_nz^n$ be given and suppose that $0 < R(f) < \infty.$ Let $m_n = o(n/\ln n),\,m_n\leq m_{n+1}\leq m_n+1,\, n\to\infty.$ Set, as before, $\pi_n:=\pi_{n,m_n}.$ Suppose now that  a subsequence $\{\pi_{n_k}\},\,n_k\in\Lambda\subset\NN$ converges $\sigma-$almost uniformly inside some  domain $U$ such that $U\supset D_{R(f)}$ and $\partial U\bigcap \partial{D_{R(f}} \not\equiv \emptyset.$ Following the classical terminology related to power series (\cite{ostrowski}), we say that  $\{\pi_{n_k}\}_{n_k\in\Lambda}$ is {\it overconvergent}.} The original  definitions and results, given for overconvergent sequences of Taylor polynomials, may be found in \cite{ostrowski}.

\smallskip \no{\bf Theorem 4, \cite{ostrowski}, {\cite{ostrowski1}}:} {\it Given a power series $f =\sum f_nz^n$  with radius of holomorphy $R_0, 0 < R_0 < \infty$ and   sequences $\{n_k\}$ and $\{n_k'\}$ with $n_k < n_k'\leq  n_{k+1},\, k = 1, 2 ...., $  suppose that

\no either

\no a) $$f_n = 0\,\, \hbox{for}\,\, n_k < n \leq n_{k'}$$ and $$n_k/n_k' \to 0,k\to\infty,\, k\to\infty.$$

\smallskip \no or

\smallskip \no b) $$\limsup n_k/n_k' < 1$$ and $$\limsup_{n\in \bigcup_{k}(n_k, n_{k'}]} |f_n|^{1/n} < 1/R_0.$$

Then

 \no a) the sequence of Taylor sums $\{S_{n_k}\}$ converges to $f$ , as $n_k\to\infty$  uniformly in the $\max-$norm inside the largest domain in $\CC$ into which $f$ is analytically continuable.%

or

\no b)  $\{S_{n_k}\}$ converges uniformly to $f$ inside neighborhoods of all regular points of $f$ on $\Gamma_{R_0}$.}

\no(here $S_n(z) = \sum_{\nu=0}^n f_\nu z^\nu.$)

 Ostrowski's theorem was extended to Fourier series associated with orthogonal polynomials in
\cite{rkk1} and to infinite series of Bessel  and of multi-index Mittag-Lefler functions in \cite{paneva}.

 \smallskip
 Before
presenting the next result, we introduce the term $G(f)$ as the {\it largest domain in $\CC$  into which $(f, D_{R_0})$ given by (1) admits a meromorphic continuation.} More exactly,  $G(f)$ is made up by the analytic continuation of the element $(f, D_{R_O})$ plus the points which are poles of the corresponding analytic function. Obviously, $D_{R(f)}\subseteq G(f).$ Further, we say that the point $z_0\in \partial D_m$ resp. $z_0\in \partial D_{R(F)}$ is {\it regular}, if $f$ is either holomorphic, or meromorphic in  a neighborhood  of $z_0.$

\smallskip \no{\bf Theorem 5, \cite{guillermo}:} {\it Let $f(z)$ be a power series with positive radius of convergence and  $m\in\NN$ be a fixed number.   Suppose that $R_m < \infty.$ Suppose that there are infinite sequences $\{n_k\}$ and $\{n_k'\}, n_k < n_k'\leq  n_{k+1}, k = 1,2,...$ such that $$\pi_{n,m} = \pi_{n_k, m}\,\, \hbox{ for}\,\, n_k < n \leq n_{k'}$$

Suppose, further,  that either

\no a) $$ \lim_k\frac{n_{k}}{n_k'} = 0 ,\, k\to\infty.$$
 \smallskip

 \no or

\smallskip

\no b)  $$\limsup_k\frac{n_{k}}{n_k'} < 1 ,\, k\to\infty.$$

Then

\no a) The sequence $\{\pi_{n_k,m}\}$ converges  to $f$, as $n_k\to\infty$,  $\sigma-$ almost uniformly inside $G(f);$

 or

\no b) $\{\pi_{n_k,m}\}$
 converges to $f$, as $n_k\to\infty$, $\sigma-$ almost uniformly in a neighborhood of each point $z_0\in \Gamma_{R(f)}$ at which $f$ is
 regular }.

 \smallskip}

{ The results of \cite{guillermo} have been extended in \cite{spain} to the $m-$th row of a large class of multipoint Pad{\'e} approximants, } associated with regular compact sets $E$ in $\CC$ and regular Borel measures supported by $E.$

 \medskip

 \no {\bf 2. Statement of the new results}

%The basic result of the present work is the following theorem, %which reflects the preceding assertions.

In the present paper, we prove

 \smallskip

\no{\bf Theorem 6:} {\it Given a power series $f$ with $R(f)\in (0,\infty)$ and a sequence of integers $\{m_n\}$ such that $m_n\leq n, m_n\leq m_{n+1}\leq m_n+1,\, m_n = o(n),\, n\to\infty$ assume that the subsequence  $\{\pi_n\}_{n\in\Lambda},\, \Lambda\subset \NN$ converges to a holomorphic, resp., meromorphic function in $\sigma-$content
  inside some domain $W$ such that
$W\bigcap  D_{R_m}^c \not\equiv \emptyset$.

  Then $$\limsup_{n\in\Lambda} \vert a_n\vert^{1/n} < 1/R(f)$$}

\no{\bf Remark:} If $m_n = m$ for all $n\in\NN$, then under the conditions of Theorem 6
 $$\limsup_{n\in\Lambda}\vert A_{n-1}\vert^{1/n} < 1/R(f).$$

\bigskip \no{\bf Theorem 7:} {\it Given the power series $f$ with  $0 < R(f) < \infty$ and a sequence of integers $\{m_n\},\, m_n\leq n, m_n\leq m_{n+1}\leq m_n+1, m_n = o(n/\log n)$ as $n\to\infty$, suppose, that $f$ is regular  at  the  point $z_0\in \Gamma_{R(f)}$. Suppose, also that there exist increasing sequences $\{n_k\}$ and $\{n_k'\},$ $n_k < n_k'\leq  n_{k+1},$ such that $\limsup_{k\to\infty}\frac{n_{k}}{n_k'} < 1$ and $\limsup_{n\in \bigcup_{k\to\infty}[n_k,n_{k}']} |a_{n}|^{1/n} < 1/R(f).$ Let  \beq\liminf_k\frac{n_{k}}{n_k'} > 0.\label{end}\eeq Then there is a neighborhood  $U$ of $z_0$ such that the sequence $\{\pi_{n_k'}\}$ converges to the function $f$ $\sigma-$almost uniformly inside  $D_{R(f)}\bigcup U.$}

\bigskip

The next result extends Theorem 5 to closed to row sequences of classical Pad{\'e} approximants.

\no{\bf Theorem 8:} {\it Let $f$ be given by (1), $0 < R(f) < \infty$ and $\{m_n\}$ be as in Theorem 7. Assume that
$ n_k < n_k' \leq n_{k+1},\, k = 1, 2, ...$ and
 \beq \pi_{n} = \pi_{n_k}\,\hbox{as}\,\,n\in\bigcup_k(n_k, n_{k}'].\label{Th02}\eeq
  Assume, further,  that either

 a)  \beq \hspace{1.2cm}n_k/n_{k'}\to 0,\,\hbox{as} \, k\to\infty\label{Th301}\eeq

\no or

b) \beq \hspace{1.2cm}\limsup_{k\to\infty} n_k/n_{k'} < 1,\,\hbox{as} \, k\to\infty\label{R}\eeq
and  $f$ is regular  at  the  point $z_0\in \Gamma_{R(f)}$.

Then

\no a)  the sequence $\{\pi_{n_k}\}$ converges to $f$ $\sigma-$almost uniformly inside $G(f)$

or

\no b) there exists a neighborhood  $U$ of $z_0$ such that the sequence $\{\pi_{n_k}\}$ converges to the function $f$ $\sigma-$almost uniformly inside  $D_{R(f)}\bigcup U.$

}

\smallskip At the end, we provide a result dealing with  overconvergent subsequences of the
 $m$th row of the classical Pad{\'e} table.

\smallskip
 \no{\bf Theorem 9:}  {\it  Let $f$ be given,
  $m\in\NN$ be fixed and  $ R_m(f):= R_m \in (0, \infty).$ Suppose that the  subsequence $\{\pi_{n_k,m}\},\, m-$ fixed,
 converges, as $n_k\to\infty$,  $\sigma-$alsmost uniformly inside a domain $U\supset D_{R_m}, \,\partial U\bigcap \Gamma_{R_m} \not\equiv
 0.$

 Then
 there exists a sequence
 $\{l_k\},\, l_k\in\NN,\, 0\leq l_k< n_k$
 such that for $n_k-l_k\leq \nu\leq n_k$ $$\limsup_{\nu\in
 {\bigcup_{k=1}^\infty} [n_k-l_k,n_k]}\vert a_{\nu}\vert^{1/\nu} < 1/R_m. $$ }

 \smallskip

\no{\bf 3. Proofs}

\medskip \no{\bf Auxiliary}

Given an open set $B$ in $\CC,$ we denote by ${\cal A}(B)$ the class of analytic and single valued functions in $B.$
 We recall that a function $g$ is meromorphic at
some point $z_0$, if there is a neighborhood $U$ of $z_0$ where $g$ is meromorphic, i.e. $g = \frac{G}{q}$ as $z\in U$, where $G\in{\cal A}(U),\, G(z_0)\not= 0$ and $q$ is a polynomial  with $q(z_0)= 0.$ We will use the notation $g\in{\cal M}(U).$

In the sequel, $D_{R},\,R > 0 $ stands for the open disk $\{z, |z| < R\};\, \Gamma_{R}:= \partial D_{R}$, respectively; $D_1:=D,\, \Gamma:= \partial D$.

With  the normalization (\ref{a1}) we have  \begin{equation}\Vert Q_n\Vert_K:= \max_{z\in K}|Q_n(z)|\leq C^{m_n},\, n\geq n_0\label{b1}\end{equation} for every compact set $K\subset \CC$, where $C=C(K)$  is independent on $n,\, 0 < C < \infty$. Under the condition $m_n = o(n),\, n\to\infty $ we have, for every $\Theta > 0$ and $n$ large enough \beq\Vert Q_n\Vert_K \leq \tilde C e^{n\Theta},\, n\geq n_0(\Theta)\label{0002}.\eeq

In what follows, we will denote by $C$  positive constants, independent on $n$ and different at different occurrences (they may depend on all other parameters that are involved) The same convention applies to $C_i, i = 1,2,....$

We take an arbitrary $\ve$ an define the open sets \beq\begin{array}{lll} \Omega_n(\varepsilon):=\bigcup_{j\leq \mu_n}(z, |z-\zeta_{n,l}|<\frac{\varepsilon}{6\mu_nn^2}),&n\geq 1\\
 \hbox{and}&\\
 \Omega(\varepsilon): =\bigcup_{n}\Omega_n(\varepsilon).\\
 \end{array}\label{begin}
 \eeq

 We have $\sigma(\Omega(\ve)) < \ve$ and $\Omega_{\ve_1} \subset \Omega_{\ve_2}$ for $\ve_1<\ve_2.$
 For any set $K\subset \CC$ we put $K(\ve):= K\setminus \Omega(\ve).$

  { Let $m_n = o(n/\ln n)$ as $n\to\infty$ and $\Theta$ be a fixed positive number. Then, as it is easy
 to check
 \begin{equation} 1/\min_{z\in K(\ve)}|Q_n(z)|  \leq   C e^{n\Theta},\,n\geq n_0(K) \label{b22} \end{equation}
 for any compact set  $K\in \CC$ and $\ve > 0$.
 If $m_n = m$ for every $n$, then
$$1/\min_{z\in K(\ve)}|Q_n(z)|  \leq   C,\, n \geq n_0(K) $$
 }

\medskip We recall in brief the  properties of the convergence in $\sigma-$content. { Let  $\Omega$ be  a domain and $\{\varphi_n\}$ a sequence of rational functions, converging uniformly in $\sigma-$content to a function $\varphi$  inside $\Omega$. If $\{\varphi_n\}\in{\cal A}(\Omega),$ then $\{\varphi_n\}$ converges uniformly in the $\max-$norm inside $\Omega.$
 If $\varphi$ has  $m$ poles in $D$, then each $\varphi_n$ has at least $m$ poles in $\Omega$;
 if each $\varphi_n$ has no more than $m$ poles in $\Omega$, then so does the function $\varphi$. For details, the reader is referred to
 \cite{Go3}}.

\medskip \no{\bf Proof of Theorem 6}

  As it follows directly from (\ref{0002})  and
from Theorem 2, \begin{equation}\limsup_{n\to\infty} \Vert P_n \Vert_K^{1/n} = 1\label{a5}\end{equation} for every compact set $K\subset D_{R(f)}$. Set $$v_n(z):= \frac{1}{n}\log{\vert\frac{ P_n(z) }{z^n}}\vert.$$ Let $\Theta$ be a fixed positive number with $e^\Theta < R(f).$ The functions $v_n$ are subharmonic in $D_{R(f)}^c$; hence, by the maximum principle (see \cite{safftotik}) and by
 (\ref{a5})
\begin{equation} v_n(z) \leq \log(\frac{{e^\Theta}}{R(f)}),\,n\in\NN,\, z\in D_{R(f)}^c, n\geq n_1. \label{a6} \end{equation}

 Let now $U_j, U_j\subset W, j = 1, 2$  be concentric open
disks of radii $0 < r_1< r_2,$ respectively, and not intersecting the closed disk $\overline{D_{R(f)}}.$

The proof will be based on the contrary to the assumption that $$\limsup_{n_k\to\infty}|a_{n_k}|^{1/n} < 1/R(f).$$
 Then there is a subsequence of $\Lambda$ which we denote again by $\Lambda$ such that
  \begin{equation}\lim_{n_k\in\Lambda}|a_{n_k}|^{1/n} = 1/R(f).\label{A}\end{equation}
Fix an $\ve, r_2-r_1 >4\ve  > 0.$ Under  the conditions of the theorem, the sequence $\{\pi_n\}_{n\in\Lambda}$ converges in $\sigma-$content inside $U_2$. Set $g$ for the limit function. Select a subsequence $\tilde\Lambda\subset \Lambda$ such that $$\sigma\{z\in U_2\setminus U_1, |\pi_{n_k'}(z) - g(z)| \geq \ve\}\leq \ve/2^k,\, n_k'\in\tilde\Lambda.$$ Set $B_{n_k'}:= \{z\in U_2\setminus U_1, |\pi_{n_k'}(z) - g(z)| > \ve\}$ and $B':=\bigcup_{n_k'\in\tilde\Lambda} B_{n_k'}.$ We have $$\sigma(\bigcup_{k+1}^{\infty}(B_{n_k'})) < \ve.$$ By the principle of the circular projection (\cite{goluzin}, p. 293, Theorem 2), there is a circle $F$, lying in  the annulus $U_{2}\setminus U_{1}$ and concentric with $\partial U_j,j=1,2$ such that $F\bigcap B' = \emptyset.$ Hence,
  $$||P_{n_k'}(z)||_F\leq C_1^{m_{n_k'}},\, n_k'\geq n_1,\,n_k'\in\tilde\Lambda$$ which yields
 $$\Vert P_{n}\Vert_{\overline U_1}\leq C_1^{m_{n}},\, n\geq n_2\geq n_1,\, n\in\tilde\Lambda$$

 Select  now a number $r$ in such
a way that the disk $D_r$ intersects the disk $U_{1}$ and set $\gamma:= \overline U_{1}\bigcap \Gamma_r.$ By construction, $\gamma$ is an analytic curve lying in the disk $U_{1}$. Applying the maximum principle to the last inequality, we get $$ \Vert P_{ n_k'}\Vert_\gamma \leq C_1^{m_{n_k'}},\,  n_k'\geq n_3\geq n_2, n_k'\in\tilde\Lambda.$$  Therefore, \begin{equation} v_{n_k'}(z) \leq  \Theta - \log r,\,\tilde n_k'\in\tilde\Lambda,\, z\in \gamma, n\geq n_4\geq n_3. \label{a7} \end{equation} Fix  now a number $\rho, R(f) < \rho$  and such that the circle $\Gamma_\rho$ does not intersect the closed disk $\overline U_2.$  By the two constants theorem (\cite{goluzin}, p. 331) applied to the domain $D_{R(f)}^c - \gamma$ there is a positive constant $\alpha =\alpha(\rho),\, \alpha < 1$ such that $$\Vert v_{n_{n_k'}} (z)\Vert_{\Gamma_{\rho}} <\alpha\Vert v_{n_k'}\Vert_{\Gamma_{R(f)}} + (1-\alpha)\Vert v_{n_k'}\Vert_\gamma,\, z\in {\Gamma_{\rho}}.$$ From  (\ref{a6}) and (\ref{a7}),  it follows that
 $$\Vert v_{n_k'}\Vert_{\Gamma_{\rho}}\leq \alpha(\log r
 - \log R(f)) + \Theta - \log r,\, n_k'\geq n_4,
n_k'\in\tilde\Lambda.$$ Hence,
 $$\limsup_{n\to\infty,n\in\tilde\Lambda}\vert v_{n}(\infty) \vert \leq \alpha(\log r
 - \log R(f)) + \Theta - \log r.$$ After letting
$\Theta$ tend to zero, we get $$\limsup_{n\to\infty, n\in\tilde\Lambda} v_{ n}(\infty) \leq \alpha(\log r
 - \log R(f))- \log r < - \log R(f).$$ The last inequality contradicts (\ref{A}), since
$$ \log\vert a_{n} \vert^{1/n}= v_{n}(\infty). $$

\no On this, Theorem 6 is proved.\hfill{Q.E.D.}

 \medskip
\no{\bf Proof of Theorem 7}

As known, the Pad{\'e} approximants are invariant under linear transformation, therefore without loss of generality, we may assume that $R(f) = 1$ and $z_0 = 1.$ Under the conditions of the theorem, there is a neighborhood of $1$, say $V$,  such that $f\in{\cal M}(V).$

Set, as before, $$\pi_{n,m_n}:= \pi_n = P_n/Q_n,$$ where $(P_n,Q_n) = 1$ and $Q_n$ are normalized as in (\ref{a1}).

Fix a positive number $\alpha > 1$ such that \beq\liminf\frac{n_k'}{n_k} > 1+\alpha,\,\, \alpha > 0.\label{Th43} \eeq
 In view of the conditions of the theorem, there is a number $\tau > 0$ such that
 \beq\limsup_{n\in \bigcup [n_k, n_k']} |a_n|^{1/n} \leq  e^{-\tau}.\label{Th44}\eeq Hence (see (\ref{c1})),  \beq \limsup_{n\in \bigcup[n_k,
 n_k'),n\to\infty}|A_n|^{1/n} \leq  e^{-\tau}. \label{a10}\eeq

Introduce the circles $C(\rho): = C_{1/2}(\rho):=\{|z-1/2| = \rho\},\, \rho > 0$ and set $ D(\rho):= \{|z-1/2| < \rho\}.$ By our previous convention,   $D_\rho:=\{z, |z| < \rho\};\, \Gamma_\rho := \partial D_\rho; D_1:=D, \Gamma_1:=\Gamma$.

Consider the function \beq\phi(R):= (\frac{1}{4R}+\frac{1}{2})^{1+\alpha}(R+\frac{1}{2})\label{0101}\eeq It is easy to  verify that there is a positive number $\delta_0$, such that $$\phi(R) < 1\, \hbox{ if}\,\, \frac{1}{2} < R < \frac{1}{2} + \delta_0.$$ Fix a number $\delta,\, 0 < \delta < \delta_0$ such that $\overline{D(\delta)}\subset D\bigcup V$ and  $\delta < e^\tau - 1.$%\label{21}\eeq With this choice of $\delta$ \beq e^{-\tau}(R + 1/2) < 1\, \hbox{for}\, R < \delta + 1/2.\label{22}\eeq

Select now a positive $\ve < \delta/4$ and introduce, as above, the sets $\Omega_n(\ve)$ and  $\Omega(\ve).$ By the principle of the circular projection, there is a number $R, 1/2 < R < 1/2 + \delta$ such that $C(R) = C(R,\ve):= C(R) \setminus \Omega(\ve).$ {Set } \beq r = \frac{1}{4R}.\label{23}\eeq

 {\sf Denote by $\omega$ the monic polynomial of
smallest degree  such that $F:=f\omega\in{\cal A}(\overline{D_{r+1/2}\bigcup D(R)}); \omega(z) = \prod_{k=1}^\mu (z - a_k),\, a_k\in \overline{D_{r+1/2}\bigcup D(R)}.$}

In what follows we will estimate the terms  $\Vert FQ_{n_k'} - \omega P_{n_k'} \Vert_{C(r)} $ and $\Vert FQ_{n_k'} - \omega P_{n_k'} \Vert_{C(R)}. $ For this purpose, we {select   a number $\Theta > 0$ such that $\Theta < \tau,\, \, e^{\Theta} (r+1/2) < 1$ and $e^{\Theta - \tau} (R+1/2) < 1.$}

By the maximum principle $$\Vert (Q_{n} F-\omega P_{n})\Vert_{C(r)} \leq \Vert Q_{n}F(z)-\omega(z)P_{n}(z)\Vert_{\Gamma_{1/2+r}},\,\leq C_0e^{n\Theta}(r+1/2)^n,\, n\geq n_0.$$ We obtain from Theorem 2, after keeping in mind (\ref{c1}), (\ref{Th43}) and the choice of $\Theta,$
 $$ \Vert Q_{n_k'}F - \omega P_{n_k'}\Vert_{\Gamma_{1/2+r}}
 \leq C_1 e^{n_k'\Theta} (r+1/2)^{n_k'}\leq C_2(e^\Theta(r+1/2))^{n_k(1+\alpha)},\, n_k\geq n_1.$$
  Thus, \begin{equation} \Vert FQ_{n_k'}-\omega P_{n_k'}\Vert_{C(r)}
  \leq C_2(e^\Theta(r+1/2))^{n_k(1+\alpha)},\,n_k\geq n_1.\label{266}\end{equation}

Estimate now $\Vert FQ_{n_k'}- \tilde\omega P_{n_k'}\Vert_{C(R)}.$

{ Clearly,  \beq\Vert F - \omega\pi_{n_k'}\Vert_{C(R)} \leq \Vert F - \omega\pi_{n_k}\Vert_{C(R)} + \Vert \omega(\pi_{n_{k}'} - \pi_{n_{k}})\Vert_{C(R)}\label{25}\eeq}

From (\ref{a5}),  we have $$\Vert \omega P_{n}\Vert_{\Gamma_{R+1/2}} \leq C_3 (e^{\Theta}(R+1/2))^{n},\, n\geq n_2.$$ On the other hand, $$\Vert \omega P_{n}\Vert_{\overline D(R)} \leq \Vert \omega P_{n}\Vert_{\Gamma_{R+1/2}}. $$ Combining the latter and the former, we get
  $$\Vert \omega P_{n_k}\Vert_{C(R)}  \leq C_3
(e^{\Theta}(R+1/2))^{n_k},\, n_k\geq n_2.$$ From here, we  obtain (see (\ref{b22})
 \beq \Vert F - \omega \pi_{n_k}\Vert_{C(R)}\leq C_4
(e^{\Theta}(R+1/2))^{n_k},\, n_k\geq n_3\geq n_2\label{26}\eeq

Let now $n_l\in [n_k, n_{k'}-1].$ %the and (\ref{b22}) and having in mind that $C(R) = C(R,\ve)$,
 By (\ref{b2}), $$\Vert \omega(\pi_{n_l+1}-\pi_{n_l}) \Vert_{C(R)} \leq \Vert \omega\Vert_{\overline D(R)}\vert A_{n_l}\vert \frac {\Vert
 z\Vert_{C(R)}^{n_l+m_{n_l}+1}}{\min_{C(R)}\vert Q_{n_l}Q_{n_l+1}\vert }$$
 which leads, thanks  (\ref{a10}) and (\ref{b22}), to
 $$\Vert \omega(\pi_{n_l+1}-\pi_{n_l}) \Vert_{C(R)} \leq C_5(e^{\Theta - \tau}(R+1/2))^{n_l}, n_k \geq n_4\geq n_3 $$
 Finally, the choice of $R$ and $\Theta$ and the conditions of the theorem imply
  $$\Vert \omega( \pi_{n_k'} - \pi_{n_k})\Vert_{C(R)} \leq \Vert\sum_{l=n_k}^{n_k'-1}\Vert \omega( \pi_{n_l+1} - \pi_{n_l}) \Vert_{C(R)}$$
  $$\leq C_6e^{n_k(\Theta - \tau)}(R+1/2)^{n_k},\, n_k\geq n_5\geq n_4 $$
  From the last inequality, combined with   (\ref{25})  and  (\ref{26}), we derive
      $$\Vert F(z)-\omega \pi_{n_k'}(z)\Vert_{C(R)}\leq C_6 (e^\Theta (R+1/2))^{n_k},\, n_k\geq n_5\geq n_4.$$ Hence, after utilization
    (\ref{0002}), we get  \begin{equation}\Vert FQ_{n_k'}- \tilde \omega P_{n_k'}\Vert_{C(R)}\leq C_7(e^\Theta)^{2n_k'}(R+1/2)^{n_k},\, n_k\geq
    n_6\geq n_5.\label{27}\end{equation}

 {We now apply  Hadamard's three circles theorem (\cite{goluzin}, p. 333, pp. 337 --
 348)
  to  $\frac{1}{n_k'}\log\vert FQ_{n_k'}(z) -\omega P_{n_k'}(z)\vert_{C(1/2)}$ and  the annulus  $\{z, r \leq |z-1/2| \leq R\}$.
  Recall that by our convention  $Rr = 1/4$.
 Using now (\ref{266}), (\ref{27}), (\ref{0101}),
 we get
  $$\frac{\log{\frac{R}{r}}}{\log\frac{1/2}{r}}\frac{1}{n_k'}\log\Vert FQ_{n_k'}- \omega P_{n_k'} \Vert_{C(1/2)}\leq
  \frac{n_k}{n_k'}(\Theta + \log \phi(R)) + 2\Theta,
 \, n_k\geq n_7\geq n_6. $$ Hence, $$\limsup_{n_k'\to\infty}\frac{1}{n_k'}\log\Vert {FQ_{n_k'}- \omega P_{n_k'}}\Vert_{C(1/2)}
  \leq \Theta(\frac{n_k}{n_k'}+2) + \frac{n_k}{n_k'}\log \phi(R)).$$ Viewing  (\ref{end}) we get,  after letting $\Theta\to 0$
   $$\limsup_{n_k'\to\infty}\frac{1}{n_k'}\log\Vert {FQ_{n_k'}- \omega P_{n_k'}}\Vert_{C(1/2)} <  0.$$
    The last inequality is strong. Hence, we
   may choose a number
  $\rho, 1/2 < \rho < R$ and close enough to $1/2$ such that  the  inequality preserves the sign; in other words, there are numbers $\rho\in
  (1/2, R)$ and $q = q(\rho) < 1$ such that $$\frac{1}{n_k'}\log\Vert FQ_{n_k'}- \omega P_{n_k'}\Vert_{\overline D(\rho)} \leq \log q, \,n_k'\geq
  n_0.$$
 From here, the $\sigma-$ almost uniform convergence inside the disk $D(R)$ immediately follows (see \cite{Go3}, Eq. (23).)}
 Indeed, fix an appropriate number $\rho.$ In view of the last inequality, $$\Vert FQ_{n_k'} - \omega P_{n_k'}\Vert_{\overline D(\rho)} \leq C'q^{n_k'},\,\, n_k'\geq n_0.$$  Take
 $\ve
 < 1/4(\rho - 1/2)$ and introduce the sets $\Omega_{n_k'}(\ve), \,n_k'> n_0$ with $\Omega_{n_0}(\ve)$ covering the zeros of the polynomial $\omega$
 (see (\ref{begin}).
 As shown above,
$\sigma (\Omega(\ve) < \ve;$ thus $ \Vert f - \pi_{n_k'} P_{n_k'}\Vert_{K(\ve)} \leq C^{"}q^{n_k'},\, n_k,\geq N. $

 On this, the $\sigma-$almost uniform convergence is established and Theorem 7 is proved. \hfill{\bf Q.E.D.}

\medskip

\no{\bf Proof of Theorem 8}

As  in the previous proof, we suppose that $R(f) = 1.$  With this convention,  \beq\limsup_{n\to\infty} \vert A_n\vert^{1/n} = 1.\label{s}\eeq

Fix a compact set $K\subset G(f)$. Our purpose is to show that $\pi_{n_k}$ converges, as $n_k\to\infty \,\sigma$-almost unformloy on $K.$ We exclude the case $K\subset D.$ In the further considerations, we assume,  that $K\nsupseteq D.$ Apparently, the generality will be  not lost.

Take a curve $\gamma_1$ such that $\gamma_1\bigcap D \not=\emptyset$, the compact set $K$ lies in the interior $B_1$ of $\gamma_1$ and $\gamma_1\subset G(f).$ Suppose that $f\in{\cal A}(\gamma_1) $ and denote by $Q$ the monic  polynomial whit zeros  at the poles of $f$ in $B_1 $ (poles are counted with their multiplicities). Set $F:= fQ;\,\hbox{deg}Q:= \mu$.   Choose a disk $B_2,\, B_2\subset D\bigcap B_1$ and not intersecting $K; \gamma_2:=\partial B_2.$ In what follows, we will be estimating $\Vert FQ_{n_k}-QP_{n_k} \Vert_{\gamma_1}$ and  $\Vert FQ_{n_k}-QP_{n_k} \Vert_{\gamma_2}.$

{ Take a  number $r_2 < 1$ such that $B_2\subset D_{r_1}$.

{  Fix $\Theta > 0$ such that $r_2e^\Theta < 1$} } Then, for every $n$ great enough there holds  $$ \Vert FQ_{n}-QP_{n} \Vert_{\gamma_2} \leq \Vert FQ_{n}-QP_{n} \Vert_{\Gamma_{r_2}} \leq C_1 (e^{\Theta} r_2)^{n+m_{n}+1},\, n > n_1$$ Hence, by (\ref{Th02}) and the choice of $r_2$ and $\Theta$  \beq \Vert FQ_{n_k}-QP_{n_k}\Vert_{\gamma_2}  =  \Vert FQ_{n_k'}-QP_{n_k'} \Vert_{\gamma_2} \leq C_1(e^{\Theta} r_2)^{n_k'}\label{th31}, n_k\geq n_1\eeq

In order to estimate $\Vert FQ_{n_k}-QP_{n_k} \Vert_{\gamma_1},$ we proceed as follows: fix a number
 $\ve, 0 < \ve < \hbox{dist}(\gamma_1,\partial G(f))/4$  %\min{\hbox{dist}(B_2,\Gamma_1)(\hbox{dist}(\gamma_1, \partial G))/4$
 and take $r_1 > 1$ such that the circle $\Gamma_{r_1}$ does not intersect the set $\Omega(\ve)$ and surrounds $\gamma_1.$

Relying on (\ref{b2}), on (\ref{b22})  and  (\ref{s}),   we get $$ \Vert \pi_{n_k} \Vert_{\gamma_1} = \Vert \sum_{n=0}^{n_k-1}\frac{A_{n}z^{n+m_n+1}}{Q_nQ_{n+1}} \Vert_{\Gamma_{r_1}}
  \leq C_3 (r_1e^{\Theta})^{n_k},\, n_k\geq n_2$$
{ Using now ( \ref{b1}) and following  the same argumentation  as in the proof of Theorem 7, we obtain } \beq\Vert FQ_{n_k} - Q\pi_{n_k}\Vert_{\gamma_1}\leq C_4(e^{\Theta}r_1)^{n_k},\, n_k\geq n_2\label{Th32}\eeq

The application of Hadamard's two constants theorem leads to $$\frac{1}{n_k}\log \Vert FQ_{n_k} - QP_{n_k}\Vert_K  \leq \alpha\frac{n_{k'}}{n_k}(\Theta + \log r_2) + (1 - \alpha)(\Theta + r_1))$$ with $\alpha:=\alpha(K) < 1.$

We get, thanks  the choice of $\Theta$  $$\lim_{n_k\to\infty}\Vert FQ_{n_k} -  Q\pi_{n_k}\Vert_K^{1/n_k} = 0$$

The statement of the theorem follows now after using standard arguments. {  On this, the proof of the first part of Theorem 8   is completed. }

\smallskip

b) The proof of the second part   { is based on the arguments provided in  } the proof of Theorem 7. As in  Theorem 7, we introduce the number $\alpha$ (\ref{Th43}), the  function $\phi(R)$ (\ref{0101}),  the circles $C(r) $ and $C(R)$ (\ref{23}) and  the polynomial $\omega$. Let $R$ and $r$ be a in THeorem 7 and set  and $F:= f\omega.$

Fix  a positive number $\Theta$ such that $e^\Theta(r+1/2) < 1, \Theta < -\frac{1}{2}\frac{\log\phi(R)}{2+\alpha}.$

We get, first, thanks (\ref{Th02}) $$\Vert FQ_{n_k} - \omega P_{n_k}\Vert_{C(r)} \leq C_1(e^\Theta(r+1/2))^{n_k(1+\alpha)},\,n_k\geq n_1$$ and, then, following the same way of considerations, $$ \Vert FQ_{n_k} - \omega P_{n_k}\Vert_{C(R)} \leq C_2(e^\Theta(R+1/2))^{n_k}, n_k\geq n_2\geq n_1.$$ Applying the tree circles theorem, we get $$\frac{\log{\frac{R}{r}}}{\log\frac{1/2}{r}} \frac{1}{n_k}\log \Vert Q_{n_k} - \omega P_{n_k}\Vert_{C(1/2)} \leq (2+\alpha)\Theta + \log \phi (R)$$ By the choice of $\Theta$, $$\frac{1}{n_k}\log \Vert Q_{n_k} - \omega P_{n_k}\Vert_{C(1/2)} < 0,\, n_k\geq n_3\geq n_2.$$ In what follows, we use standard arguments to complete the proof of (b), Theorem 8. \hfill{\bf Q.E.D.}

\smallskip

\no{\bf Proof of Theorem 9}

 Without losing the generality, we
assume that $R_m = 1$ and $\tau_n = 0$ for all $n.$
 Normalize the polynomials $Q_n$ as
it was done in (\ref{a1}) with $R(f)$ replaced by $R_m$. Fix a positive number $\ve,\, \ve < 1/2$ and introduce the set $\Omega(\ve).$ Select a number $R,\, R > 1$ such that $\Gamma_R\bigcap \Omega(\ve) = \emptyset$. Recall that (see (\ref{b2})) there are positive constants $C_j(\ve):= C_j, j = 1,2$ such that \beq \frac{n^{2m}}{C_1} \leq \min_{z\in\Gamma_R}|Q_n(z)| < \Vert Q_n\Vert_{\Gamma_R} \leq C_2,\, n \geq n_0.\label{new3}\eeq In the sequel, we assume that $C_1,\,C_2 > 1.$

By Theorem 1, there is  a positive number   $\tau = \tau(R)>0$ such that
 \begin{equation}\Vert P_{n_k}\Vert_{\Gamma_R}\leq C(e^{-\tau}R)^{n_k},\, n_k\geq n_1\geq n_0;\label{new1} \end{equation} and (by the maximum
 principle for subharmonic functions),  \begin{equation}\vert a_{n_k}\vert\leq C(e^{-\tau})^{n_k}\label{new2},\, n_k\geq n_1\geq n_0 \end{equation}
 Without losing the generality, we suppose that \beq R^{m+1}\leq C\leq C_1.\label{000}\eeq

 We will prove that for every $l, 0 \leq l \leq n_k$ and for $n_k$ great enough \beq\Vert P_{n_k-l}\Vert_{\Gamma_R} \leq (2C_2)^{l}
 C_1^{l+1}(e^{-\tau}R)^{n_k}\prod_{j=0}^{l-1}(n_k-j)^{2m},\, n_k \geq n_2.\label{new00}\eeq

From the last  inequality, it follows directly that
 \beq|a_{n_k-l}| \leq (2C_2)^{l} C_1^{l+1}R^l(e^{-\tau})^{n_k}\prod_{j=0}^{l-1}(n_k-j)^{2m}, n_k \geq n_2\label{001}\eeq

 We prove first (\ref{000}) for $l = 1.$ For this purpose, we introduce the polynomial $${\cal P}_{n_k}:= P_{n_k-1}Q_{n_k} - P_{n_k}Q_{n_k-1}.$$ By
 definitions of Pad{\'e} approximants (see (\ref{b2})), \beq {\cal P}_{n_k}(z) =
 A_{n_k}z^{n_k+m+1}\eeq where, according to (\ref{c1}), $$A_{n_k-1} = a_{n_k}\prod(\frac{-1}{\tilde\zeta_{n_k-1, l}})$$ (recall that by presumption
 the  defect $\tau_n = 0$ for all $n.$)  Viewing  (\ref{new1}) and (\ref{000}), we get
  \beq  \Vert {\cal P}_{n_k}\Vert_{\Gamma_R}\leq |a_{n_k}|{R}^{n_k+m+1}\leq C_1(e^{-\tau}{R})^{n_k},n_k\geq n_3\geq n_1,\,\,\tau:=
 \tau(R).\label{new5}\eeq
 Keeping now
 track of (\ref{new1}) and (\ref{new3}), we arrive at $$ \Vert P_{n_k-1}Q_{n_k} \Vert_{\Gamma_R} \leq C_1(e^{-\tau}{R})^{n_k} +
 C_1C_2(e^{-\tau}{R})^{n_k},\, $$ which yields  \beq \Vert P_{n_k-1} \Vert_{\Gamma_R} \leq (C(e^{-\tau}R)^{n_k} +
 C_1C_2(e^{-\tau}R)^{n_k})/\min_{\Gamma_R} |Q_{n_k}(z)|,\, n_k \geq n_4\geq n_3\label{new4}\eeq $$\leq 2C_1C_2^2(e^{-\tau}{R})^{n_k}n_k^{2m}.$$  We
 further get \beq |a_{n_k-1}| = |\frac{1}{2\pi i}\int_{\Gamma_R}\frac{P_{n_k-1}(z)}{z^{n_k}}dz|\leq 2C_1C_2^2Re^{-\tau
 n_k}n_k^{2m},n_k\geq n_4\,\label{new6}\eeq

 Suppose now that (\ref{new00}) is true for $l-1, l\geq 2.$ In other words,
 $$\Vert P_{n_k-l+1}\Vert_{\Gamma_R} \leq (2C_2)^{l-1} C_1^{l}(e^{-\tau}R)^{n_k}\prod_{j=0}^{l-2}(n_k-j)^{2m},\, n_k\geq n_5$$
 and
 $$|a_{n_k-l+1}| \leq (2C_2)^{l-1} C_1^{l}R^l(e^{-\tau})^{n_k}\prod_{j=0}^{l-2}(n_k-j)^{2m},\,n_k\geq n_5$$

 Introducing into considerations the polynomial ${\cal P}_{n_k-l+1}:= P_{n_k-1}Q_{n_k-l+1} - P_{n_k-l+1}Q_{n_k-1}$ and following the same arguments
 as
 before, we  see that  (\ref{new00}) and (\ref{001}) are true also for $l.$

Equipped with inequality (\ref{001}), we complete the proof of the theorem. We will be looking for numbers $l$ such that $$|a_{n_k-l}|^{1/n_k-l} < 1.$$ Set $2RC_2C_1:= C4.$ We check that $$\log |a_{n_k - l}|^{1/(n_k-l)}\leq \psi_{n_k}(l), $$ where $$\psi_{n_k} (x):= \frac{C_4x + C_1}{n_k-x} + \frac{x2m\log n_k}{n_k-x} - \tau \frac{n_k}{n_k-x}.$$ For $n_k$ large enough, say $n_k\geq n_6$, $\psi$ is strongly increasing, and $\psi(0) < 0.$ Hence, there is a number $x_k \in (0, n_k)$ such that $\psi_{n_k} (x) < \psi_{n_k} (x_k) < -\tau/2$ every time when $0 < x < x_k.$ Set $l_k:= x_k.$ Therefore  $$\limsup_{{n\in\bigcup_{k=1}^\infty}[n_k-l_k, n_k]}|a_{n}|^{1/n} < 1$$  \hfill{\sf Q.E.D.}

%{\bf According to previous considerations, Theorems 6 - 9 will be valid for mulripoint PA of Newton type.}

 %newpage


\begin{thebibliography}{20}
  \bibitem{gonchar1} A. A. Gonchar, Poles of the rows of the  Pad{\'e} table
 and meromorphic continuation of functions. Mat. Sb. 115(157)(1981), 590 -- 615. English transl. in Math. USSR Sb. 43(1982).
 \bibitem{Pe} O. Perron, Die Lehre von den Kettenbr{\"u}chen, {
Teubner, Leipzig}, 1929.
  \bibitem{Go3} { A. A. Gon{ch}ar, On the convergence of generalized Pad{\'e} approximants of meromorphic functions, { Mat. Sbornik,} { 98} (140)
    (1975), 564 -- 577, { English translation in Math. USSR Sbornik,} { 27} (1975), No. 4, 503--514.}
    \bibitem{montessus} R. de Montessus de Ballore, Sur le fractions continues algebriques, Bull. Soc. Math. France 30(1902), 28 -- 36.
    \bibitem{blkov} H. P. Blatt, R. K. Kovacheva, Growth behavior and zero distribution of rational functions, Constructive Approximation, 34(3),
        (2011),     393 -- 420.
    \bibitem{blatt1} H. P. Blatt, Convergence in capacity of rational approximants of meromorphic functions, Publ. Inst. Math. (N.S.) 96 (110),
        (2014), 31 -- 39.
\bibitem{blatt} H. P. Blatt, Overconvergence of Rational Approximants of Meromorphic Functions, Progress in Approximation Theory and Applicable
    Analysis,   Springer - Verlag, 2017.
    \bibitem{goluzin} G.M. Goluzin, Geometric Theory of Function of Complex Variables, in Russian, Publ. House Nauka, Moskow, 1966.

    \bibitem{baker} G. A. Baker, Jr., P. Gr. Morris,   Pad{\'e} Approximants, Second Edition, Encyclopedia of Mathematics and Applicatiobs, V. 59,
        Cambridge University  Press, 1996.
     \bibitem{vavprsue} V. V. Vavilov, G. Lopes Lagomasino,  V. A. Prokhorov, On an inverse problem for the rows of the Pad{\'e} table. Mat. Sb.
         110(152)(1979),
        117 -- 129;  English transl.   in  Math. USSR. 38(1981).
    \bibitem{ostrowski} A. Ostrowski, {\"U}ber Potenzreihen, die {\"u}berkonvergente Abschnittsfolgen besitzen. Deutscher Mathematischer Verein, 3
    2(1923), 185 -- 194.
    \bibitem{ostrowski1} A. Ostrowski, {\"U}ber eine Eigenschaft gewisser Potenzreihen mit  unendlich vielen verschwindenden Koeffizienten,
        Berliner         Berichte 1921, 557 -- 565.
  \bibitem{guillermo} A. Fern{\'a}ndes Infante, G. Lopez Lagomasino, Overconvergence of subsequences of
  rows of Pad{\'e} approximants with gaps,
      Journal     of Computational and Applied mathematics 105(1999), 265 -- 273.
 \bibitem{spain}    B de la Calle Ysern, JM Ceniceros,  Rate of convergence of row sequences of multipoint Pad{\'e} approximants,
  Journal of Computational and Applied Mathematics, 284(2015), 155 -- 170.
       \bibitem{safftotik} E. B. Saff, V. Totik, Logarithmic Potentials with External Fields, Springer-Verlag, Berlin, Heidelberg, 1997.
        \bibitem{landkoff} N. S. Landkov, Foundations of Modern Potential Theory, Grundlehren der mathematischen Analysis, Springer-Verlag,
            Berlin, Heidelberg, 1972.
        \bibitem{rkk1} R. K. Kovacheva, Zeros of sequences of partial sums and overconvergence, Serdica Mathematical Journal, vol. 34, 2008, 467
            -- 482.
\bibitem{paneva}  Jordanka Paneva-Konovska, From Bessel to Multi-index Mittag-Leffler Functions: Enumerable Families, Series in Them and
    Convergence, World Scientific Publishing Europe,  2016.

\end{thebibliography}
\end{document}